\documentclass[12pt,leqno]{article}
\usepackage{latexsym,amsfonts,epsfig}
\input epsf
\newtheorem{theorem}{Theorem}[section]

\newtheorem{lemma}[theorem]{Lemma}

\begin{document}
\title{When does third order efficiency imply fourth order efficiency}

\author{Shanti Venetiaan\thanks{Institute of Graduate Studies and Research, Anton de Kom Universiteit van Suriname, Paramaribo, Suriname} }

\maketitle

\begin{abstract}
In this article third and fourth order efficiency are studied in the framework of translation equivariant location estimators. We assume $X_1,...,X_n$ i.i.d. $f(\cdot -\theta)$. By recognizing that equality in a special form of the Cauchy-Schwarz inequality leads to a certain dependence of the cumulants of the maximum likelihood estimator (MLE) for $\theta$, it is shown that this MLE is fourth order efficient if the underlying distribution is Gumbel. Contrary to similar results which were previously published this result is not based on symmetry.
\end{abstract}

t


\pagestyle{myheadings} \thispagestyle{plain} \markboth{...}{{\small THIRD AND FOURTH ORDER EFFICIENCY}}

\section{Introduction}
Let $X_1,...,X_n$ be independent and identically distributed random variables with common distribution function $F(\cdot - \theta)$, $\theta \in \mathbb R$, and with density $f(\cdot -\theta)$, with respect to Lebesgue measure on $(\mathbb R,B)$. We assume that the derivative $f'(\cdot - \theta)$ is such that the Fisher information for location is finite, i.e.
\begin{equation}\label{eqno1}
I(f) = \int (\frac{f'}{f})f < \infty
\end{equation}
holds.
We estimate the location parameter $\theta$ by a translation equivariant estimator $T_n = t_n(X_1, ...,X_n)$, i.e. $T_n$ satisfies 
\begin{equation}
t_n(x_1+a,..., x_n+a) = t_n(x_1,..., x_n) + a,\quad a,x_1,...,x_n \in \mathbb R\nonumber
\end{equation}
We denote the distribution function of $a_nT_n$, under $f(\cdot)$ by $G_n(\cdot)$, namely
\begin{equation}\label{Hn}
G_n(y) = P_{\theta}(a_n(T_n-\theta)\leq y),\quad y\in \mathbb R,
\end{equation}
and we are interested in constructing confidence intervals for $\theta$ based on $T_n$. The confidence intervals are of the form
\begin{equation}\label{confint}
[T_n - \frac{1}{a_n}G_n^{-1}(v), T_n - \frac{1}{a_n}G_n^{-1}(u)] 
\end{equation}
with $0 < u < v < 1$.
It may be shown that $G_n$ has a density (c.f. Theorem 1.1. of Klaassen (1984)). From this follows that $G_n(G_n^{-1}(u)) =u$ for all $u$, $0<u<1$, and hence intervals like (\ref{confint}) have coverage probability $v-u$. The length of the confidence interval in (\ref{confint}) is $\frac{1}{a_n}(G_n^{-1}(v)-G_n^{-1}(u))$ and one may be interested in using the length of a confidence interval as a measure of the performance of the estimator $T_n$. Shorter intervals with the same coverage probability will be considered to be linked to better estimators than estimators related to longer intervals. The ideas mentioned above have already been considered in Klaassen and Venetiaan (2010). In that article the confidence interval inequality was introduced which says that under certain regularity conditions,

\begin{equation}
G_n^{-1}(v) - G_n^{-1}(u) \ge \epsilon\label{ineq}
\end{equation}
for all $n \in \mathbb N$, $u$, $v$, $0<u<v<1$, and $c >0$ and $\epsilon > 0$ satisfying the relations
\begin{equation}
u = P_0 \left(\prod_{i=1}^n f\left(X_i +
\frac{\epsilon}{a_n}\right) \ge c\, \prod_{i=1}^n f(X_i)\right)
\end{equation}
and
\begin{equation}
v =P_{-\epsilon/a_n} \left(\prod_{i=1}^n f\left(X_i +
\frac{\epsilon}{a_n}\right) \ge c\, \prod_{i=1}^n f(X_i)\right).
\end{equation}

Asymptotic expansions were derived for both the left hand side and the right hand side of (\ref{ineq}) and it can be shown that the inequality in (\ref{ineq}) may be used to prove that first order efficiency implies second order efficiency. This was first proved by Pfanzagl (1979) and has after that been studied by scholars like Pfanzagl and Wefelmeyer(1985), and  Bickel, Chibisov and van Zwet (1981). In Klaassen and Venetiaan (1994) the spread inequality of Klaassen(1984) was used to prove this phenomenon. For the maximum likelihood estimator which is first order efficient and therefore second order efficient,  Klaassen and Venetiaan (2010) also found that locally even third order efficiency is obtained by the MLE and this third order efficiency immediately implies fourth order efficiency, because of the symmetry of the intervals and the even polynomials in the fourth order term. Ghosh (1994) introduced the conjecture that third order efficiency implies fourth order efficiency and this assertion was proved by Akahira(1996) by studying the coverage probability of symmetric confidence intervals.

In this paper we will limit ourselves to the case of taking $T_n$ to be the MLE and $a_n = \sqrt{nI(f)}$ in (\ref{Hn}).
In the proof that locally third order efficiency implies fourth order efficiency by Klaassen and Venetiaan(2010) the difference $G_n^{-1}(v) - G_n^{-1}(u)-\epsilon$ is studied for a special choice for $u$ and $v$. If this special choice is not made, studying the difference gives that third order efficiency and even fourth order efficiency may be obtained when a certain  dependence among the cumulants of the MLE holds. This dependence of the cumulants may be viewed as equality in a certain case of the Cauchy-Schwarz inequality, and results in
 
\begin{equation}
\frac{f'}{f}(X) = \lambda\left(\frac{f''}{f}(X)-(\frac{f'}{f}(X))^2+I(f)\right)
\end{equation}
for some $\lambda \in \mathbb R$.
Note that $\frac{f''}{f}(X)-(\frac{f'}{f}(X))^2$ is the derivative of $\frac{f'}{f}(X)$, so we have a differential equation and its solution is the Gumbel distribution with density $f(x) = \frac1{\beta}\exp[{\frac{x-\alpha}{\beta}-\exp{\frac{x-\alpha}{\beta}}}]$, for some $\alpha \in \mathbb R$, and $\beta >0$.
We will also show that when the underlying distribution is Gumbel, the MLE estimator for location is fourth order efficient. In other words, regardless of the choice of $u$ and $v$, an MLE estimator for location is fourth order efficient if the underlying distribution is Gumbel.

In Section 2 we will repeat the theorem on the confidence interval inequality, the theorem concerning the Cornish Fisher expansion for $G_n^{-1}(u)$ and also repeat the theorem which gives an expansion for the $\epsilon$ in (\ref{ineq}). The main result is stated in Section 3 and its proof may be found in the same section.

\section{Confidence interval inequality and asymptotic expansions}
\setcounter{equation}{0}
In this section we mention some results which were published in Klaassen and Venetiaan (2010) and Venetiaan (2010).
\begin{theorem}\label{stelling1}
Let $X_1, ..., X_n$ be independent and identically distributed
random variables with common distribution function $F(\cdot -
\theta)$, $\theta \in R$, and with density $ f(\cdot - \theta)$,
with respect to Lebesgue measure on $(\mathbb R, {\cal B})$. Let
$T_n= t_n(X_1,...,X_n)$ be a translation equivariant estimator for
$\theta$ and let $H_n$ be the distribution function of $a_nT_n$,
$a_n > 0$, under $f(\cdot)$, i.e.
\begin{equation}\label{eq:no4}
H_n(y) =  P_0 \left( a_n T_n\le y \right) = P_\theta \left(a_n(T_n
- \theta)\le y \right) ,\quad y \in \mathbb R\,.
\end{equation}
Fix $u$ and $v$ with $0<u<v<1$, and assume that there exist
$\epsilon >0$ and $c >0$ satisfying the relations
\begin{equation}\label{epsilonc1}
u = P_0 \left(\prod_{i=1}^n f\left(X_i +
\frac{\epsilon}{a_n}\right) \ge c\, \prod_{i=1}^n f(X_i)\right)
\end{equation}
and
\begin{equation}\label{epsilonc2}
v =P_{-\epsilon/a_n} \left(\prod_{i=1}^n f\left(X_i +
\frac{\epsilon}{a_n}\right) \ge c\, \prod_{i=1}^n f(X_i)\right).
\end{equation}
If $H_n^{-1}(\cdot)$ is strictly increasing at $u$ and continuous
at $v$, then
\begin{equation}\label{epsilon}
H_n^{-1}(v) - H_n^{-1}(u) \ge \epsilon.
\end{equation}
\hfill $\Box$
\end{theorem}

Furthermore we define our notion of asymptotic efficiency. Asymptotic efficiency of $T_n$ means that $H_n$ converges weakly to the standard normal distribution function $\Phi$ as $n \to \infty$,
\begin{equation}
H_n^{-1}(u) = \inf \{y \in R: H_n(y) \ge u \} \to \Phi^{-1}(u),\quad 0 < u < 1.
\end{equation}
By efficiency of the $j$-th order we mean that the expansions up to and including the $n^{-(j-1)/2}$ term have optimal coefficients. As mentioned earlier, we take $T_n$ to be the maximum likelihood estimator, $a_n = \sqrt{nI(f)}$ in (\ref{Hn}) which choice leads to asymptotic efficiency of $T_n$. Furthermore in our framework, optimal coefficients for the $n^{-(j-1)/2}$ term of $G_n^{-1}(\cdot)$ will be obtained when  $|G_n^{-1}(v) - G_n^{-1}(u)-\epsilon| = o(n^{-(j-1)/2})$.

We now introduce some notation.
\begin{eqnarray}\label{etas}
&& \eta_2 = E\psi_2^2(X_1)/I^2(f),\,\, \eta_3= E\psi_1^3(X_1)/I^{3/2}(f),\,\, \eta_4=E\psi_1^4(X_1)/I^2(f)\,,\nonumber\\
&& \eta_5 = E\psi_1^5(X_1)/I^{5/2}(f)\,,\,\, \mbox{and}\ \
\eta_6= E\left(\psi_2(X_1)\psi_3(X_1)\right)/I^{5/2}(f)
\end{eqnarray}
are defined with
\begin{equation}
\psi_i(x) = \frac{f^{(i)}(x)}{f(x)}\,, \quad x \in \mathbb R,
\end{equation}
The next two theorems were previously published in Venetiaan (2010) and Klaassen and Venetiaan (2010) respectively.
The first one gives an asymptotic expansion for $G_n^{-1}(\cdot)$ which may be used for the left hand side of (\ref{ineq}).

\begin{theorem}\label{stelling2}
Let $X, X_1,...,X_n$ be i.i.d. with common density $f(\cdot -
\theta_0)$. Let $\rho(\cdot) = - \log f(\cdot)$ satisfy the following conditions.

\begin{enumerate}
\item
For all $K \subset \mathbf{R}$ compact, $\sup_{\theta\in K}
E_{\theta} \rho^2(X) = A < \infty$.
\item
$\rho(\cdot)$ is five times differentiable.
\item
There exists a function $R(\cdot)$ and a $\delta > 0$ such that,
for every $y \in \mathbb R$, $|\theta| < \delta $:
\[ |\rho^{(5)}(y) -\rho^{(5)}(y-\theta)|\le R(y)|\theta|\quad\mbox{and}\quad E_0 R^{5/2}(X) < \infty\,.\]
\item
$E_0 |\rho^{(i)}(X)|^5 < \infty$ for $i = 1,...,5$.
\end{enumerate}
Then 
\begin{equation}
|G_n^{-1}(v)-\tilde G_n^{-1}(v)| = o\left(\frac{1}{n{\sqrt n}}\right)\label{thmgn}
\end{equation}
holds with
\begin{eqnarray}\label{eq:no10}
\tilde G_n^{-1}(v) &=& z_v\left[ 1 + \frac 1{12\sqrt n}\,\eta_3 z_v \right.\nonumber\\
& & \qquad \quad +
\frac1{72n}\left\{(-9+12\eta_2-\eta_3^2-5\eta_4)z_v^2
-9-2\eta_3^2+3\eta_4\right\} \\
&& \qquad \quad+ \frac1{144n{\sqrt n}}
\left\{(6\eta_2\eta_3-3\eta_3-\frac{19}{12}\eta_3^3-\eta_3\eta_4+
\frac{24}5 \eta_5-18\eta_6)z_v^3 \right. \nonumber\\
&&\left. \left. \qquad \qquad \qquad \qquad
+(12\eta_2\eta_3-15\eta_3-\frac{67}9\eta_3^3
+3\eta_3\eta_4-\frac95\eta_5)z_v\right\}\right] \nonumber
\end{eqnarray}
 and $z_v = \Phi^{-1}(v)$. In other words the inverse of the distribution function of $T_n$ admits the Cornish-Fisher expansion 
$\tilde G_n^{-1}(\cdot)$.
\hfill $\Box$
\end{theorem}

The next theorem results in an asymptotic expansion for $\epsilon$ in (\ref{ineq}).

\begin{theorem}\label{stelling3}
Let the conditions of Theorem \ref{stelling1} and \ref{stelling2} be satisfied. Then for
\begin{eqnarray} \tilde\epsilon& =& z_v - z_u +\frac{\eta_3}{12\sqrt{n}}(z_v^2 - z_u^2)\nonumber\\
&\quad+ &\frac1{288n}\left[ (12\eta_2+5\eta_3^2-8\eta_4)(z_v^3-z_u^3)\right.\nonumber\\
&\qquad +&(36 -36\eta_2+9\eta_3^2+12\eta_4)(z_v^2z_u-z_u^2z_v)\nonumber\\
&\qquad + &\left.(-36-8\eta_3^2+12\eta_4)(z_v -z_u)\right]\label{epstilde}\\
&\quad + &\frac1{12960n^{3/2}}\left[
(270\eta_2\eta_3+60\eta_3^3-180\eta_3\eta_4+162\eta_5-540\eta_6)
(z_v^4 - z_u^4)\right.\nonumber\\
&\qquad + &(-540\eta_2\eta_3+540\eta_3+270\eta_3^3-270\eta_5+1080\eta_6)(z_v^3z_u- z_u^3z_v)\nonumber\\
&\qquad +&\left.
(-270\eta_3-400\eta_3^3+630\eta_3\eta_4-162\eta_5)
(z_v^2- z_u^2)\right],\nonumber
\end{eqnarray}
we have 
\begin{equation}\label{eps}
|\epsilon -\tilde \epsilon| =o(1/n\sqrt n),
\end{equation}where $z_u = \Phi^{-1}(u)$

\hfill $\Box$
\end{theorem}

\section{Main result and proof}
In this section we state the main result of this paper. 
\begin{theorem}\label{stelling4}
Let the conditions of Theorem \ref{stelling1} and \ref{stelling2} hold. Then,
\begin{equation}
|G_n^{-1}(v) - G_n^{-1}(u) -\epsilon| = o(\frac1{n\sqrt n}),\quad 0<u<v<1,
\end{equation}
 if and only if the underlying distribution is Gumbel or $v= 1-u$.
In other words, in the framework of Theorem \ref{stelling1} and Theorem \ref{stelling2} the maximum likelihood estimator for location is efficient up to fourth order if and only if the underlying distribution is Gumbel or if $v= 1-u$.
\hfill $\Box$
\end{theorem}

In the proof of Theorem \ref{stelling4} we wil study the difference $\tilde G_n^{-1}(v) - \tilde G_n^{-1}(u) -\tilde\epsilon$ 
which equals
\begin{equation}
-\frac1{96n}(12-12\eta_2+3\eta_3^2+4\eta_4)[z_v^3-z_u^3+z_uz_v^2-z_u^2z_v] + o(\frac1{n})
\end{equation}
It's clear that optimal third order coefficients are obtained for $\tilde G_n^{-1}(v) - \tilde G_n^{-1}(u)$, i.e. third order efficiency may be obtained when $12-12\eta_2+3\eta_3^2+4\eta_4 = 0$. When we divide both sides of this expression by $3$ and take a closer look, we see that this expression is a case of equality in the Cauchy-Schwarz inequality. We may state the following lemma.

\begin{lemma}\label{lemma}
Use the framework of the theorems above. Then
\begin{equation}
\eta_3^2= 4\eta_2 -\frac43\eta_4-4\label{CSeta}
\end{equation}
if and only if 
the underlying distribution is Gumbel.
 \hfill $\Box$
\end{lemma}

\subsubsection*{Proof of Lemma \ref{lemma}.}

We first mention that 

\begin{eqnarray}
E\psi_1\psi_2(X_1) &=& \frac12 E\psi_1^3(X_1),\quad E\psi_1^2\psi_2(X_1) = \frac23 E\psi_1^4(X_1),\nonumber\\
E\psi_1(X_1) &=& 0, \quad \mbox{and}\quad E\psi_2(X_1) =0.
\end{eqnarray}

Note that the following form of the Cauchy-Schwarz inequality:

\begin{equation}
(EYZ)^2 \le EY^2EZ^2\label{CS}
\end{equation}
with
$Y= \psi_1$ and $Z= (\psi_1'-E\psi_1')$ results in
\begin{eqnarray}\label{remark}
[E\psi_1(\psi_1'-E\psi_1')]^2 &\le& E\psi_1^2E(\psi_1'-E\psi_1')^2\nonumber\\
{[E(\psi_1(\psi_2 -\psi_1^2+I))]}^2 &\le& IE(\psi_2 - \psi_1^2+I)^2\nonumber\\
{[E(\psi_1\psi_2 -\psi_1^3+I\psi_1)]}^2 &\le& IE(\psi_2^2 - 2\psi_1^2\psi_2+2I\psi_2+ \psi_1^4-2I\psi_1^2+I^2)\nonumber\\
\frac{[E\psi_1^3]^2}{4}&\le & IE(\psi_2^2-\frac{\psi_1^4(X)}{3}-I^2)\nonumber\\
\frac{\eta_3^2}{4}&\le& \eta_2 -\frac{\eta_4}{3}-1,
\end{eqnarray}
where $\psi_1$, $\psi_1'$ and $I$ are short for $\psi_1(X_1)$, $\psi_1'(X_1)$, and $I(f)$ respectively. Note that (\ref{CSeta}) corresponds to equality in the Cauchy Schwarz inequality. It's well known that equality in the Cauchy-Schwarz inequality (\ref{CS}) is obtained if and only if $Y =\lambda Z$ for some $\lambda \in \mathbb R$. For our case this means that equality in (\ref{remark}) is obtained if and only if $\psi_1 = \lambda((\psi_1'-E\psi_1'))$
 holds. Now note that this is a differential equation and it's solution is the density of the Gumbel distribution.
\hfill $\Box$(end of proof of Lemma \ref{lemma}).

\subsection{Proof of Theorem \ref{stelling4}}
Assume that $|G_n^{-1}(v) - G_n^{-1}(u) - \epsilon| = o(1/n\sqrt n)$, then we get
\begin{eqnarray}
|\tilde G_n^{-1}(v) - \tilde G_n^{-1}(u) -\tilde \epsilon|&\le& |\tilde G_n^{-1}(v) -G_n^{-1}(v)| + |G_n^{-1}(u) -\tilde G_n^{-1}(u)|\nonumber\\
&\quad +& |\epsilon -\tilde \epsilon| + |G_n^{-1}(v) - G_n^{-1}(u) -\epsilon|\nonumber\\
&\le& o(\frac1{n\sqrt n})
\end{eqnarray}
by using (\ref{thmgn}), (\ref{eps}), and the assumption. As $\tilde G_n^{-1}(v)-\tilde G_n^{-1}(u)-\tilde\epsilon$ only has terms up to and including a $1/n\sqrt n$-term, this means that $\tilde G_n^{-1}(v)-\tilde G_n^{-1}(u)-\tilde\epsilon = 0$. In other words, we have $\eta_3^2= 4\eta_2 -\frac43\eta_4-4$ or we have $z_u = -z_v$. As Lemma \ref{lemma} states $\eta_3^2= 4\eta_2 -\frac43\eta_4-4$ only holds when the underlying distribution is Gumbel. The other case $z_u = -z_v$ means that $v = 1-u$.  

Now, on the other hand, let $v = 1-u$, then $z_u = -z_v$ and it is easy to see that the $1/n$-term and $1/n\sqrt n$-term of $\tilde G_n^{-1}(v) - \tilde G_n^{-1}(u) -\tilde\epsilon$ vanish.  But then, 

\begin{eqnarray}\label{triangular}
|G_n^{-1}(v) - G_n^{-1}(u)-\epsilon| &\le& |G_n^{-1}(v) -\tilde G_n^{-1}(v)|+ |G_n^{-1}(u)-\tilde G_n^{-1}(u)|\nonumber\\
&\quad+&  |\epsilon-\tilde \epsilon|+|\tilde G_n^{-1}(v) - \tilde G_n^{-1}(u)-\tilde\epsilon |\nonumber\\
&=& o(\frac1{n\sqrt n})
\end{eqnarray}
because of (\ref{thmgn}), (\ref{eps}).

At last, we assume that the underlying distribution is Gumbel and we see that the Fisher information and the relevant cumulants are 

\begin{eqnarray}
&& I(f) =1,\,\,\eta_2 = 5,\,\, \eta_3= -2,\,\, \eta_4=9\,,\nonumber\\
&& \eta_5 = -44\,,\,\, \mbox{and}\ \
\eta_6= -13,
\end{eqnarray}
Note that substitution of these quantities in $\tilde G_n^{-1}(v) - \tilde G_n^{-1}(u) -\tilde \epsilon$ results in $0$. In the same fashion as in (\ref{triangular}) it may be shown that $|G_n^{-1}(v) - G_n^{-1}(u)-\epsilon|= o(1/n\sqrt n)$.
\hfill $\Box$ (end of proof of Theorem \ref{stelling4}).

\bigskip
\noindent{\bf Acknowledgments.}
The author thanks Chris Klaassen for pointing out that $\eta_3^2/4 \le \eta_2-\eta_4/3 -1$ is a special case of the Cauchy Schwarz inequality and for his feedback with regards to this paper.

\noindent 

\bigskip

\noindent{\bf\Large References}

\begin{description}

\item
Akahira, M. (1996).
\newblock Third order efficiency implies fourth order efficiency:
a resolution of the conjecture of J.K. Ghosh.
\newblock Annals of the Institute of Statistical Mathematics 48:365--380.



\item
Bickel, P.J. Chibisov, D.M. van Zwet, W.R. (1981).
\newblock On Efficiency of first and second order.
\newblock International Statistical Review 49:169--175.







\item
Ghosh J. (1994).
\newblock {\em Higher Order Asymptotics}.
\newblock NSF-CBMS regional conference series Probability and Statistics, Institute of Mathematical Statistics.



\item
Klaassen, C.A.J. (1984).
\newblock Location estimators and spread.
\newblock The Annals of Statistics 12:311--321.


\item
Klaassen C.A.J. Venetiaan, S.A. (1994).
\newblock Spread inequality and efficiency of first and second order.
\newblock In: Hu\u skov\'a, M. Mandl, P. ed.,
\newblock {\em Asymptotic Statistics, Proceedings of the Fifth Prague Symposium}:341--348.
\newblock Heidelberg: Physica-Verlag.

\item
Klaassen C.A.J. Venetiaan, S.A. (2010).
\newblock Optimizing lengths of confidence intervals: Fourth-order efficiency implies in location models.
\newblock Communications in Statistics-Theory and methods 39, 8-9:1437--1448


\item
Pfanzagl, J. (1979).
\newblock First order efficiency implies second order
efficiency.
\newblock In: Jure\u ckov\'a, J. {\em Contributions to
Statistics: J. H\'ajek Memorial Volume}:157--196.
\newblock Prague: Academia.

\item
Pfanzagl, J. Wefelmeyer, W. (1985).
\newblock {\em Asymptotic Expansions for General Statistical Models,
Lecture Notes in Statistics} 31.
\newblock New York: Springer Verlag.


\item
Venetiaan, S.A. (2010).
\newblock An expansion for the maximum likelihood estimator of location and its distribution function. 
\newblock Brazilian Journal of Probability and Statistics 23,1: 82--91.


\end{description}
\enddocument